\documentclass[%
 reprint,
 amsmath,amssymb,
 aps,
]{revtex4-1}
\usepackage{graphicx}
\usepackage{bm}
\usepackage{amssymb}
\usepackage{amsfonts}
\usepackage{amsmath}
\usepackage{color}
\usepackage{graphicx}
\usepackage{overpic}
\usepackage{epstopdf}
\usepackage[font=footnotesize,labelfont=bf]{caption}
\usepackage[labelformat=simple,font=footnotesize,labelfont=bf]{subcaption}

\usepackage{placeins}

 \usepackage{amsthm}
 \usepackage{psfrag}

\DeclareMathOperator{\tr}{tr}

\def\beq{\begin{equation}}
\def\eeq{\end{equation}}

\begin{document}


\title[Ergodicity in umbrella billiards]{Ergodicity in umbrella billiards}
\author{Maria F. Correia}
\affiliation{
CIMA-UE, Department of Mathematics, University of \'Evora, Rua Rom\~ao Ramalho, 59, 7000-671 \'Evora, Portugal.
}%
\email{mfac@uevora.pt}

\author{Christopher Cox}
\affiliation{
Department of Mathematics, Washington University in St. Louis, One Brookings Drive, St. Louis, MO 63130-4899
}%
\email{clcox@wustl.edu}
\author{Hong-Kun Zhang}%
\affiliation{Department of Mathematics and Statistics, University of Massachusetts, 710 North Pleasant Street, Amherst, MA 01003.
}
\email{hongkun@math.umass.edu}

\date{\today}

\begin{abstract}
We investigate a three-parameter family 
of billiard tables with circular arc boundaries. These umbrella-shaped billiards may be viewed as a generalization of two-parameter moon and asymmetric lemon billiards, in which the latter classes comprise
instances where the new parameter is $0$.
Like those two previously studied classes, for certain parameters umbrella billiards exhibit evidence of chaotic behavior despite failing to meet certain criteria for defocusing or dispersing, the two most well understood mechanisms for generating ergodicity and hyperbolicity.  For some parameters 
corresponding to non-ergodic lemon and moon billiards, small 
increases in the new parameter transform elliptic $2$-periodic points into a cascade of higher order elliptic
points. These may either stabilize or dissipate as the new parameter is increased. We characterize the periodic points and present evidence of new ergodic examples.

\end{abstract}

\maketitle

\section{Introduction}
Billiards are dynamical systems which are based on a simple model but which nonetheless provide deep physical insights and pose fundamental questions in statistical mechanics, quantum mechanics, and broadly across many branches of physics. On a bounded region  $Q\subset \mathbb{R}^2$ (the billiard table), an infinitesimal particle moves along segments at unit speed, changing direction according to the law of specular reflection upon collisions at boundaries.
The essential link in billiards between the geometry of the table and the dynamics of the system facilitates a robust model which has proved useful in approaching problems ranging from the foundations
of the Boltzmann's ergodic hypothesis [\onlinecite{Bu74}], to the description of shell effects in semiclassical physics [\onlinecite{BB97}], to the design of microwave resonators in quantum chaos [\onlinecite{J99}], and many other other applications [\onlinecite{AFMKR08,DOP,Gri09, HeTo,HW}].
In particular, ergodic properties are determined by the shape of the table, producing a spectrum of behaviors from completely integrable to strongly chaotic. Many questions remain unanswered in the presence of non-integrable dynamics, and it is upon such questions that we will concentrate. 

The origins of the field of chaotic billiards may be traced to Sinai [\onlinecite{S70}], who established the ergodicity of {\em dispersing} billiards and opened the door to many previously unapproachable problems. The {\em defocusing mechanism} of Bunimovich [\onlinecite{Bu74}] extended the field of study to convex tables including the well-known {\em stadium}, which he demonstrated to be hyperbolic and ergodic; and {\em flower} billiards, closely related to the billiard classes of interest in this paper. Further elucidation of defocusing billiards  followed from  Wojtkowski [\onlinecite{Woj86}], Markarian [\onlinecite{Ma88}],
Donnay [\onlinecite{Do91}] and Bunimovich [\onlinecite{Bu92}].  More recently, Bunimovich and Grigo  [\onlinecite{BG10}]  conjectured that absolute focusing is a necessary requirement for a typical convex table to be ergodic.

We are particularly interested in billiards which are not dispersing and  do not meet any known defocusing criterion, but which nonetheless exhibit chaotic properties.  Among the known examples most relevant to the new class of billiards we will investigate are {\em annular} billiards, introduced by Saito et al.\ [\onlinecite{Saito82}] and later extensively studied in [\onlinecite{AFMKR08,Boh93,Gou01,HR02,DF16}]. Benettin and Strelcyn [\onlinecite{BS78}] looked at one-parameter {\em oval} tables and observed notable properties including bifurcation phenomenon, the coexistence of elliptic and chaotic regions, and the separation of the chaotic region into several invariant components. In [\onlinecite{DRW96}] the ovals were generalized to a two-parameter family encompassing seven varieties, including special cases of {\em lemon}, {\em moon}, and a particular example of a class which in this paper we will designate as {\em umbrella} billiards, while  [\onlinecite{BHHS}] gives an alternate generalization of [\onlinecite{BS78}] to {\it squash} billiard tables, on which the elementary defocusing mechanism does not take place. 
In [\onlinecite{HeTo}] symmetric lemon billiards were considered, and
recently a class of asymmetric  lemon-shaped convex billiard tables were constructed in [\onlinecite{CMZZ13}], obtained by intersection of two disks in the plane. It was also proved that  a subclass of these billiards are indeed hyperbolic using continued fraction techniques [\onlinecite{BZZ14}]. These, along with the moon billiards recently investigated by the authors in  [\onlinecite{CZ15}], are the direct antecedents of the current investigation.

\begin{figure}[t]

\hspace{-.9 cm}
\includegraphics[width=60mm]{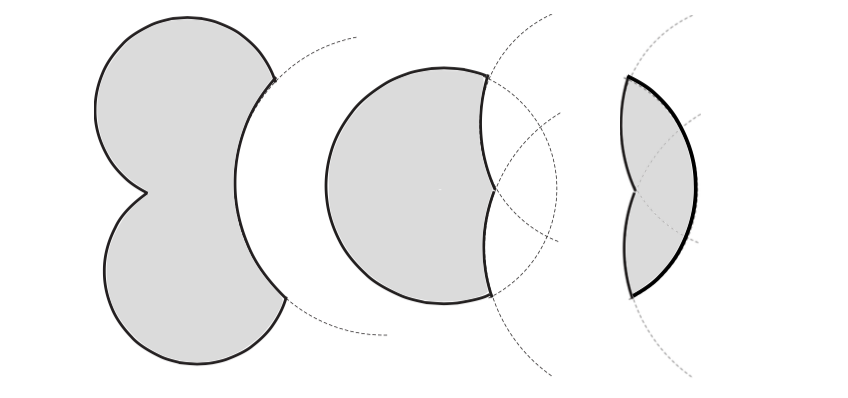}
\caption{\small{Three umbrella billiards, created by duplicating a framing disk of  a moon billiard ($Q_M^2(1.3,1.3,1.6)$, left and $Q_M^1(1,1.4,1)$, center) or a lemon billiard ($Q_L(1.3,1.5,1.2)$, right) and separating centers of the new disks.}}
\label{umbrellatypes}
\end{figure}

 Lemons and moons may be identified parametrically as $Q(B,R)$, where a circle of radius $R\geq 1$ overlaps a unit circle with centers separated by $B>0$, using the central or outer regions to form lemon or moon tables respectively. Umbrella tables  $Q_L(B,R,B_1)$ and $Q_M^i(B,R,B_1)$ are formed by duplicating the unit circle and separating the centers of the new-formed (initially overlapping) disks to a distance $B_1\geq0$ units.   (See Figure \ref{umbrellatypes}; see Section \ref{sec:umbrella} for details.) For the moon type, the superscript will distinguish between the construction in which the circle corresponding to the dispersing edge is duplicated (Type 1) and the case in which the circle corresponding to the focusing edge is duplicated (Type 2). For the lemon-based umbrellas no such distinction is needed, though the asymmetry will result in non-unique parametrizations.  
 
\begin{figure}[t]

\hspace{-1cm}
\includegraphics[width=80mm]{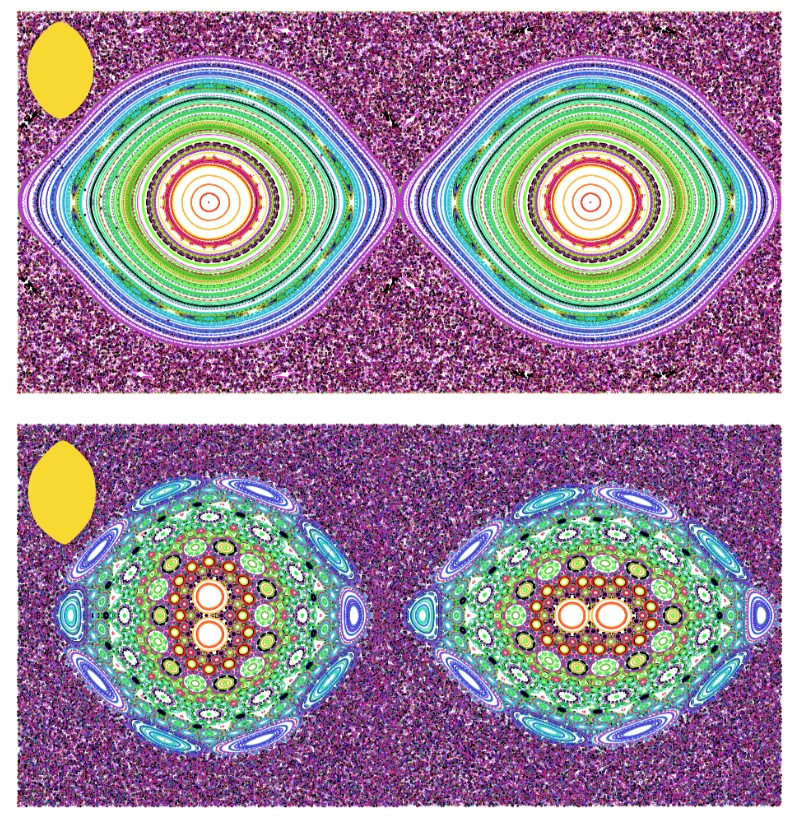}
\caption{\small{The phase portrait of the lemon billiard $Q_L(1.0.75,0)$ (above) and of the umbrella billiard $Q_L(1.0,0.75,0.05)$ (below) with outlines of the tables (upper left insets). From $B_1=0$ to $B_1=0.05$ the change in the table is minimal, but the dynamical change is significant.}}
\label{cascade}
\end{figure}

Even with this simple construction, in which all components are circular arcs, our numerical results show that these billiards still enjoy rich ergodic and chaotic properties, and small modifications may result in notable differences in the dynamics, as in Figure \ref{cascade}. Depending on the combination of both initial conditions and  parameters, the phase spaces present a rich structure
which contains invariant spanning curves, Kolmogorov-Arnold-Moser (KAM) islands and chaotic seas, and suggest that subclasses may be completely ergodic relative to the standard billiard measure. 

\begin{figure}[b]

  \includegraphics[width=60mm]{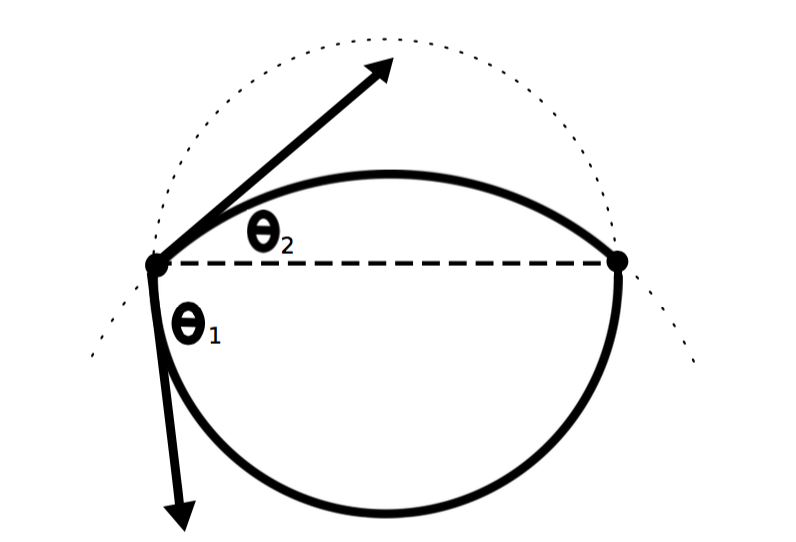}
\caption{Two-arc billiards, including lemons, moons, and two-petal flowers, may be parametrized by tangent angles of the arcs relative to the central axis.}
\label{fig:thetaparam}

\end{figure}

\begin{figure*}[t]
\hspace{-1cm}
\includegraphics[width=165mm]{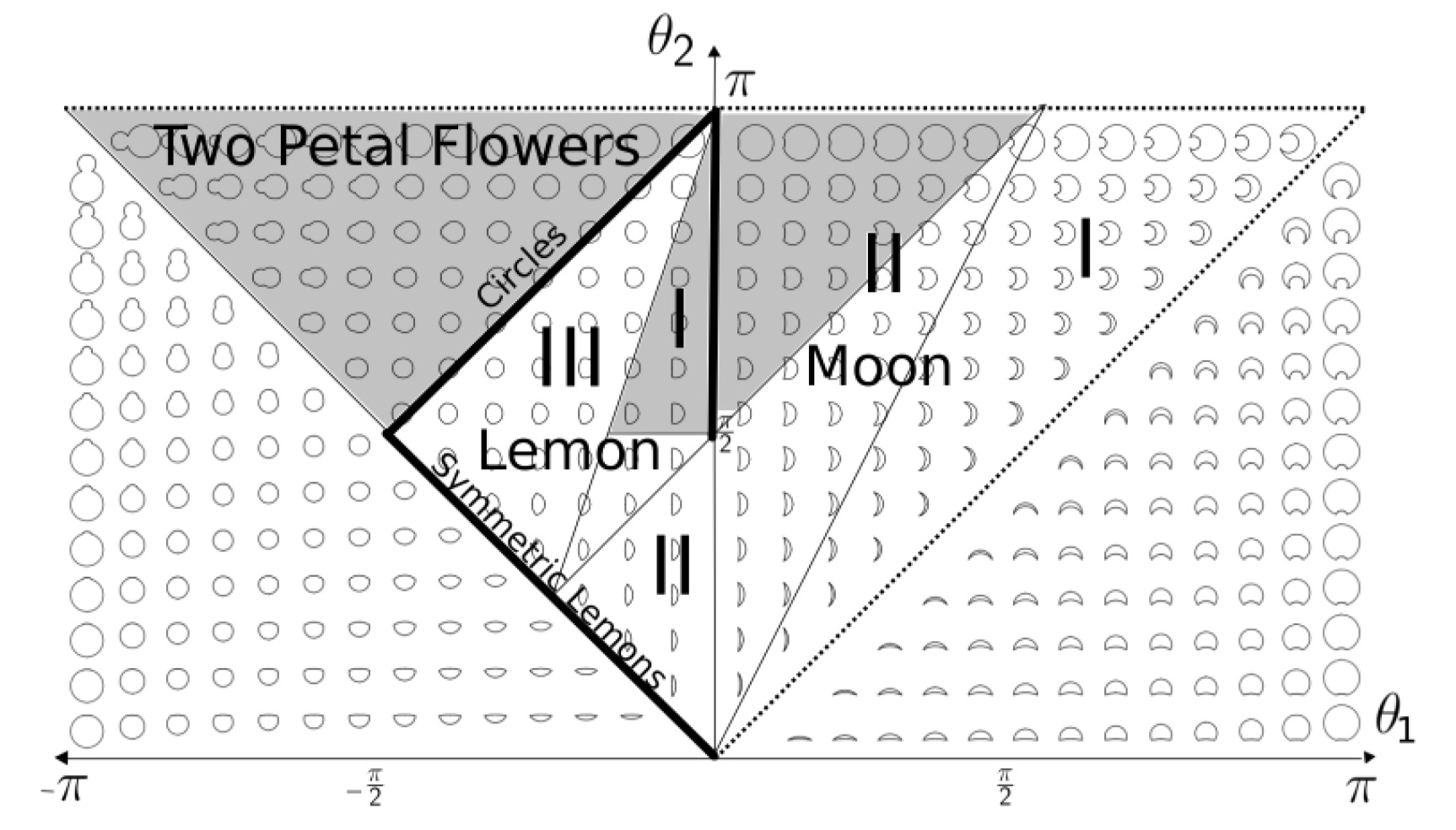}
\caption{\small{The class of billiards formed by two circular arcs can be completely parametrized in the upper triangular region. (The numbered subregions are the translations into the $\theta$ parametrization of the regions given in [\onlinecite{CMZZ13}] and [\onlinecite{CZ15}].) 
The two-petal flowers in the upper left are analytically known to be ergodic. Lemon billiards  (middle left) and moon billiards (right) exhibit numerical evidence of ergodicity in some regions, shown in gray. Under the umbrella modification, with $B_1>0$, the lemon type ergodic region expands from I into III and the moon type expands through region II and (for sufficiently large $B_1$) into region I. See Section \ref{sec:phase} for a more precise description of the moon case.}}
\label{2theta}
\end{figure*}

We introduce an alternative to the $(R,B)$ parametrization which has the advantage of uniting the moon and lemon families in the larger class of two-arc billiards. Aligning the vertices along a horizontal axis and letting $\theta_1$ and $\theta_2$ be the signed angle of the tangents of the two arcs at the left vertex (Figure \ref{fig:thetaparam}), we normalize by scaling the $\theta_1$ circle to unit radius. After reducing through identifications by relabeling and symmetry the a priori $2\pi$ square parameter space reduces to one triangular quadrant.  (Note that each $\theta$-plane corresponds to  a $B_1$, the third parameter which is unaltered, and in cases where $B_1 \neq 0$ we will use the $\theta$ values associated to the base case, not the modified angles.) The family of all billiards in this triangle includes not only asymmetric lemons and moons but also the two-petal variety of flower billiards, a class for which the ergodicity was established analytically under the defocusing mechanism of Bunimovich.  Figure \ref{2theta} summarizes the known and new billiard tables viewed through the lens of this parametrization.

In Section \ref{sec:umbrella} we discuss the antecedent classes and describe the construction of umbrella billiards, and in Section \ref{sec:periodic} we discuss periodic points, looking at the cases where the umbrella billiard mirrors their base types as well as examples where they diverge markedly. The last two sections give numeric evidence of new ergodic billiards, considering Poincar\'e surfaces of sections and the transition of periodic points from elliptic to hyperbolic in Section \ref{sec:phase}, and  Lyapunov exponents in Section \ref{sec:lyap}.

\section{Umbrella Billiards}
\label{sec:umbrella}

We begin this section with a brief description of the construction and dynamics of the parametric families of moon billiards and asymmetric lemon billiards, then introduce the generalized umbrella class. 

\begin{figure}[b]

\hspace{-1cm}

\includegraphics[width=60mm]{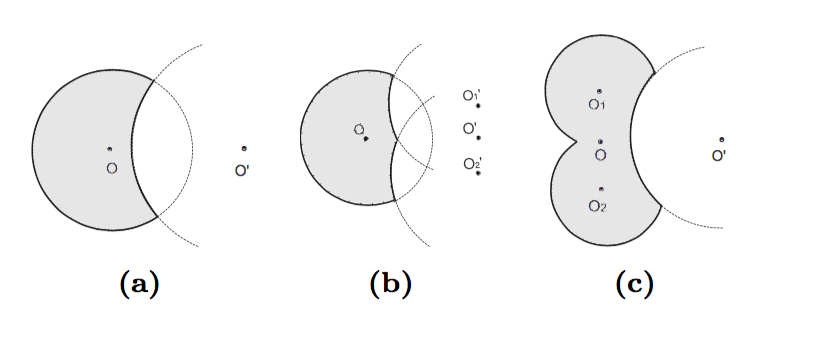}

\caption{\small{A family of umbrella billiards obtained by modifying moon billiards: (a) standard moon billiards; (b) the first variation, replacing the dispersing edge; (c) the second variation, replacing the focusing edge.}}
\label{D1D1DR}
\end{figure}

\begin{figure}[b]

\hspace{-1cm}

\includegraphics[width=60mm]{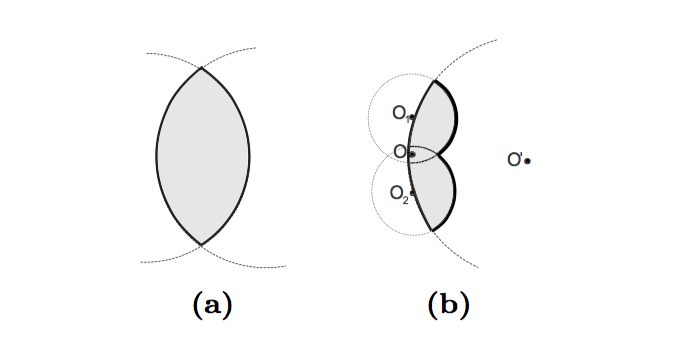}

\caption{\small{A family of umbrella billiards obtained by modifying lemon billiards: (a) standard lemon billiards; (b) the lemon-type umbrella billiard.}}
\label{perlemon}
\end{figure}

Starting with a unit disk 
and a  second disk  of radius $R\geq 1$, with the two centers separated by distance $B$, moon billiard tables are comprised of the unit disk minus the overlap with the radius $R$ disk (Figure \ref{D1D1DR}a). This defines a two-parameter family of moon-shaped billiard tables with boundary made of two circular arcs. Alternatively, an asymmetric lemon billiard table is comprised of the intersection of the two disks (Figure \ref{perlemon}a).
 We use the designation $Q_M(R,B,0)$ for moons and $Q_L(R,B,0)$ for lemons, replacing the now ambiguous $Q(R,B)$ used for the former class in [\onlinecite{CZ15}] and the latter in  [\onlinecite{CMZZ13}].

\begin{figure*}[t]

  \includegraphics[width=180mm]{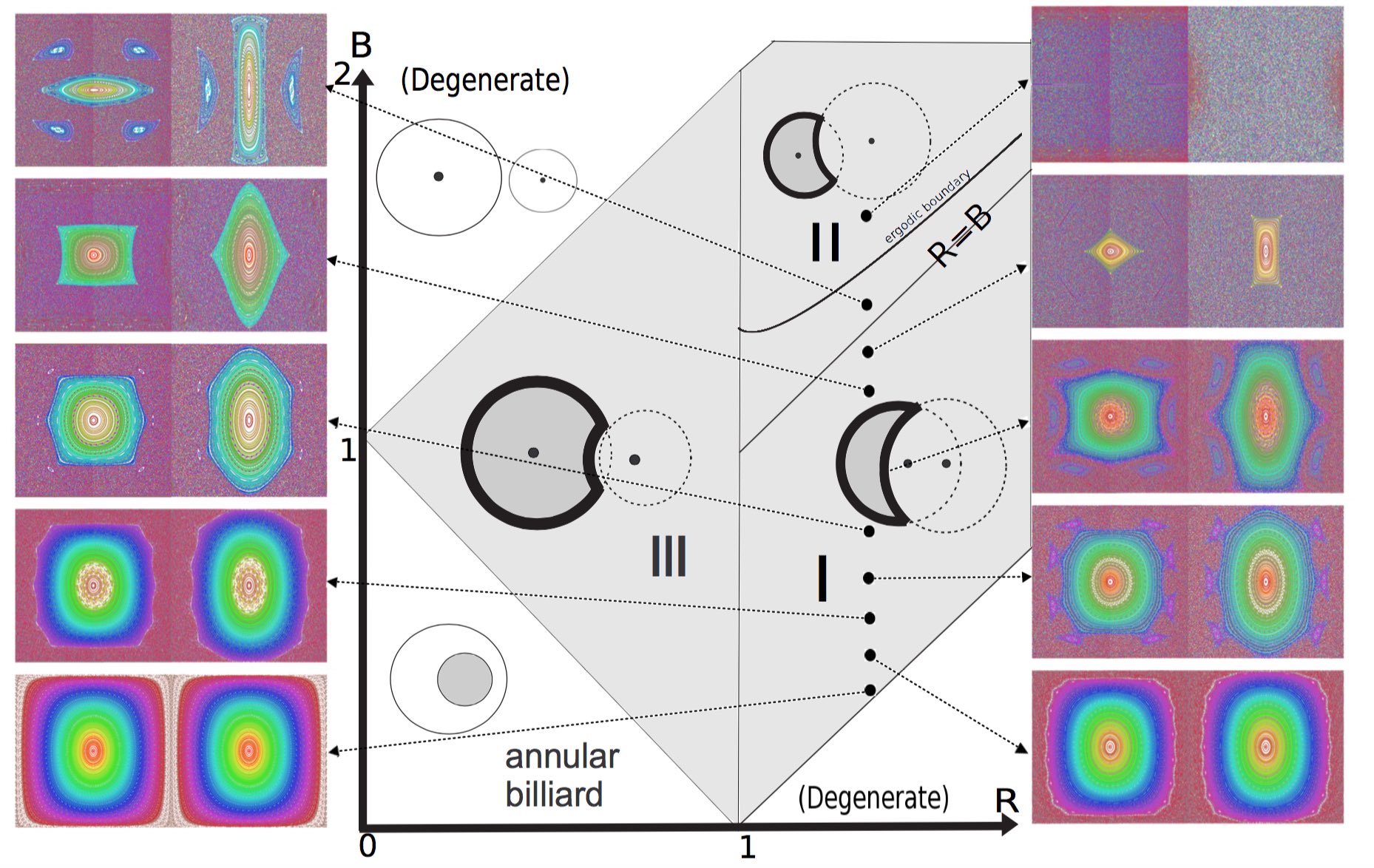}
\caption{The $R,B$ parameter space (center) for moon billiards, with phase portraits for $R=1.2$ and the indicated values of $B$ ranging from $0.21$ to $1.5$. Examples in Region II above the curve, where neither circle contains the center of the other, appear to be ergodic, while for smaller $B$ stable period two elliptic points dominate the dynamics.}
\label{moonoverview}

\end{figure*}

In both cases the shape of the boundary, and consequently the dynamics,  is controlled by the parameters $R$ and $B$. For moon billiards, typical examples have one or more elliptic points coexisting with a chaotic region [\onlinecite{CZ15}]. When $B$ is small and the centers of the two framing circles are close, elliptic islands about 2-periodic points are prominent; in contrast, sufficient separation $B$ relative to the radius $R$ results in apparently ergodic behavior. (See Region I and Region II in Figure \ref{moonoverview}.)

\FloatBarrier

\begin{figure*}[t]
  \includegraphics[width=180mm]{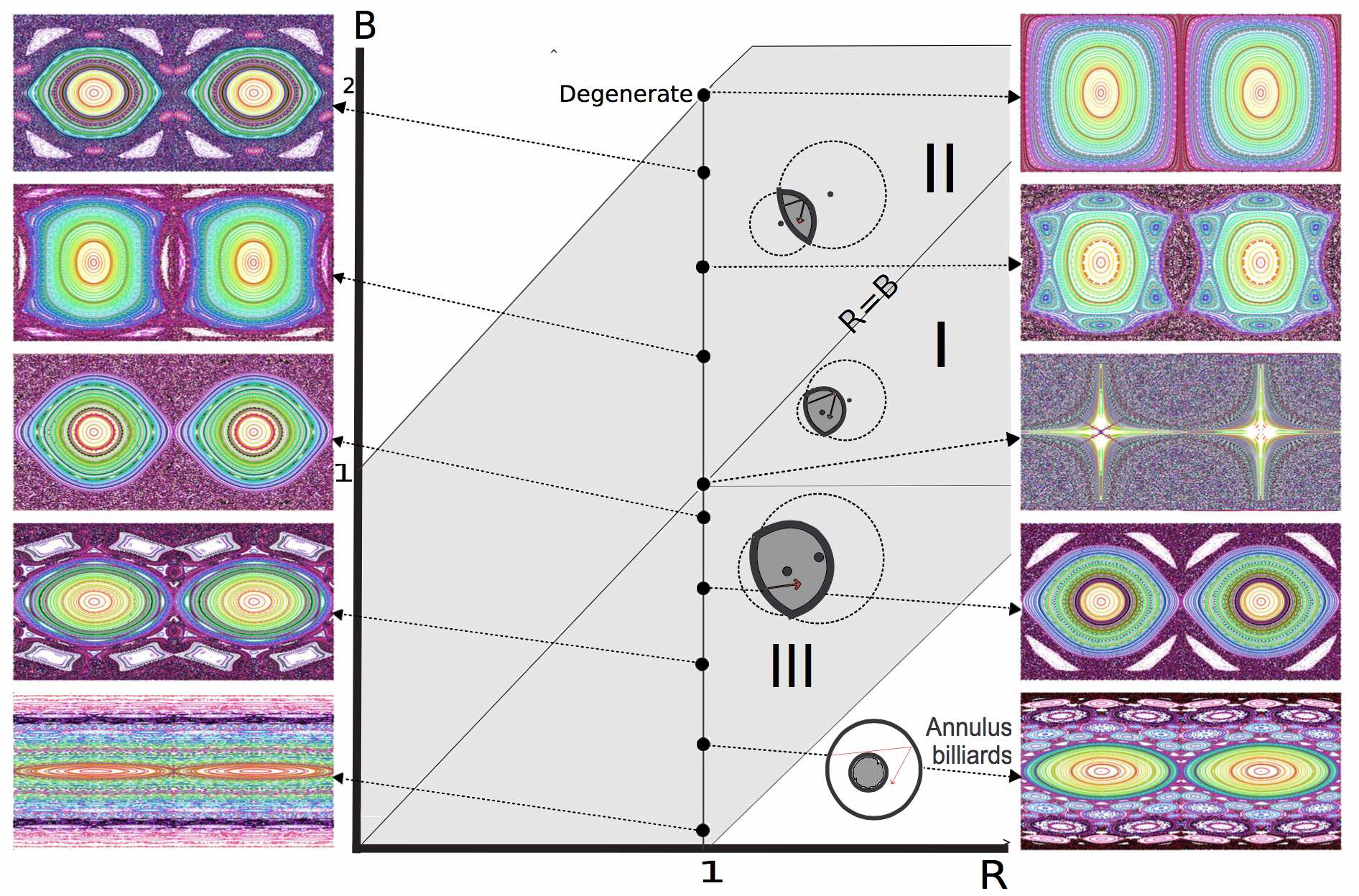}
\caption{The $R,B$ parameter space (center) for lemon billiards. For symmetric lemon tables with $R=1$ there is an isolated parabolic case $Q(1,1)$ (middle right) while all other cases are elliptic. However, the asymmetric instances are frequently hyperbolic. Region I appears to be entirely ergodic. }
\label{lemonoverview}
\end{figure*}

Symmetric lemon billiards may be parametrized by a single parameter $B$ (the subclass $Q_L(1,B,0)$ in the $RB$ parametrization) and as first observed in [\onlinecite{HeTo}] display elliptic behavior about a central $2$-periodic point for all instances except the isolated parabolic example $Q_L(1,1,0)$. As $B$ approaches $2$ and the table becomes thin, the elliptic islands around the central 2-periodic point expand and the surrounding chaotic region becomes narrow. Conversely, when $B$ approaches $0$ and the table approaches the circular case, a host of elliptic islands emerge, becoming long and narrow approaching the integrable case. If $R\neq1$, however, and asymmetric billiards are considered, hyperbolicity [\onlinecite{BZZ14}] and apparent ergodicity [\onlinecite{CMZZ13}] often arise, 
as in Region I in Figure \ref{lemonoverview}.

To create the first two modified classes, start with the moon billiard construction. Duplicate the radius $R$ disk to create two disks, and create the new table by allowing the two centers to move away symmetrically in opposite directions orthogonal to the axis between the center of the original radius $R$ disk and the unit disk. (See Figure \ref{D1D1DR}b). Let $B_1>0$ be the distance between the centers, and designate this first moon variation $Q_M^1(R,B,B_1)$.
For the second variation of the moon billiard, parallel the construction duplicating the radius $R$ disk instead, thereby altering the focusing edge of the original moon billiard (Figure \ref{D1D1DR}c). This type will be designated $Q_M^2(R,B,B_1)$.

A similar construction may be applied to obtain an umbrella variation of the lemon billiard, duplicating the unit circle and moving the two copies apart, but using the intersection of the unaltered disk with the union of the new disks instead as the table. 
See Figure \ref{perlemon}b. Denote this class by $Q_L(R,B,B_1)$.

There are three corners on these tables, which break down the smoothness of the boundary 
and will lead to the existence of nontrivial singularity curves.
More precisely, the singularity set of this table $Q(R,B)$ consists of three vertical segments in the phase space 
based at the three corner points, as well as the horizontal lines corresponding to grazing.
Additionally, we limit the class to nondegenerate modifications, specifically cases in which $B_1$ is sufficiently small that the billiard table consists of a single, simply connected region. The maximum allowable $B_1$ varies by type and by the parameters $R$ and $B$, but in all cases $B_1 < 2.$

Each type of umbrella billiard may be extended into a modification partitioning a moon or lemon billiard edge into not merely two but any number of new edges. The general $n$-umbrella billiard is by obtained by replacing the duplicated circles with an arbitrary number of circles, spaced evenly along the axis.
To round out the modifications of the class of two-arc billiards, an identical construction might be applied to two-petal flowers. It is clear that all non-degenerate cases would be three-petal flower billiards that satisfy Bunimovich's defocusing mechanism, and thus are analytically established as ergodic. Accordingly, our focus will be on investigating the ergodicity of lemon and moon types, starting with the base cases and then considering umbrella billiards.   

\FloatBarrier

\section{Characteristics of periodic orbits}
\label{sec:periodic}

In this section we consider the periodic points appearing in the base lemon and moon classes as well as the modified class of umbrella billiards.  For any fixed billiard table $Q$, let $\mathcal{M}$ be the space of unit vectors based at the boundary $\partial Q$, pointing inwards, and endowed with natural topology. $\mathcal{M}$ can be viewed as a closed cylinder with nature coordinates $x=(s, \theta)$, where $s\in [0,|\partial Q|]$ is the arc-length parameter on the wall $\partial Q$, oriented counterclockwise, and $\theta\in [0,\pi]$ is the angle formed by the vector $x$ and the positive tangent direction to $\partial Q$ at the base point $x$.
The set $\mathcal{M}$ is a natural cross-section of the phase space for the billiard flows. The first return map (or the Poincar\'e map) obtained by restricting the flow on $\mathcal{M}$ is called the discrete billiard map, $T: \mathcal{M}\to  \mathcal{M}$, $T(s,\theta)=(s_1,\theta_1)$. The billiard map preserves the probability measure $\mu$ on $\mathcal{M}$ with $d\mu=c\sin\theta\, ds\, d\theta$, where $c$ is a normalizing constant.   A periodic point $x=T^kx$ is then said to be hyperbolic, parabolic and elliptic if $|\tr(D_{x}T^k)|>2$ (unstable), $|\tr(D_{x}T^k)|=2$ (neutrally stable) and $|\tr(D_{x}T^k)|<2$ (stable) respectively. (See for instance [\onlinecite{Ber81}].)

Of particular interest, then, are types of periodic points which may transition from elliptic to hyperbolic as the parameters vary. This consideration will inform the discussion in Section \ref{sec:phase}, when we look for boundaries in the umbrella billiards parameter space where the elliptic points dissipate and the transition to apparent ergodicity occurs.
For certain parameters of all classes under discussion, 2-periodic orbits colliding within a single circular arc may occur. For moon and umbrella-moon types, such orbits may tangentially graze the boundary at one (or for $n$-umbrellas more) points. As the former is parabolic and the latter hyperbolic (see [\onlinecite{CZ15}]) they are omitted from further discussion.

\subsection{Periodic points of lemon type}

\begin{figure}[htb]
\centering
\includegraphics[width=80mm]{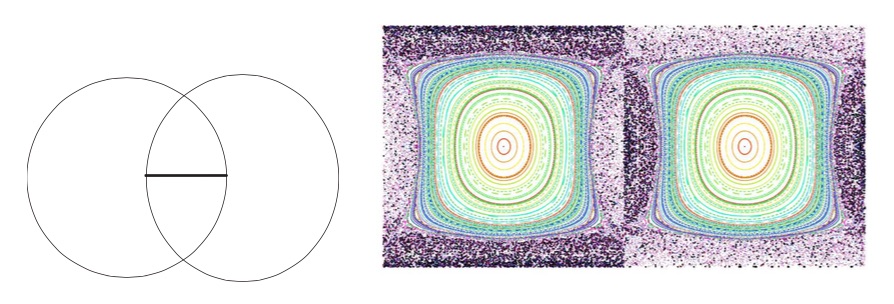}
\caption{The dominant 2-periodic orbit (left) and the phase space (right) of $Q_L(1,1.35,0)$. }
\label{d-135}
\end{figure}

For many parameters, the phase space of lemon billiards is dominated by elliptic periodic points. In Figure \ref{d-135}, a central 2-periodic point is contained in an elliptic island surrounded by a smaller chaotic sea. The orbits of the type appearing in this figure will persist for many larger values of $B$, but the shapes of the corresponding islands undergo some interesting transformations.

\begin{figure}[htb]
\centering
\includegraphics[width=80mm]{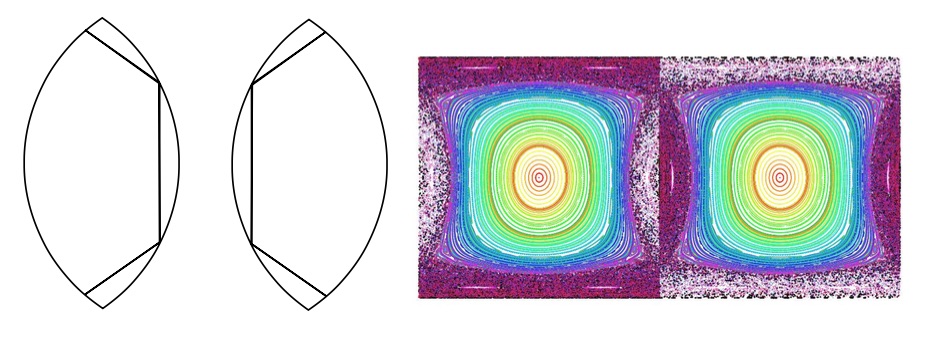}

\caption{\small{For $Q_L(1,1.37,0)$, the orbits corresponding to the outlying islands (left) and the phase space (right).}}
\label{d-137}
\end{figure}

\begin{figure}[b]
\centering
\includegraphics[width=65mm]{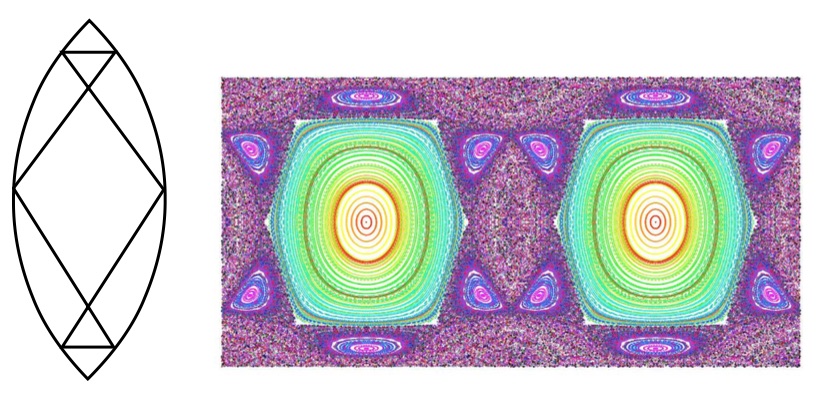}
\caption{Left: Period 6 orbit for the table of $Q_L(1,1.58,0)$. Right: Phase space of $Q_L(1,1.58,0)$.}
\label{d-158}
\end{figure}

New elliptic islands start to form when $B$ goes over  $1.35$. In Figure \ref{d-137}, we illustrate two-sided periodic trajectories around the central elliptic island island and two pairs of periodic six orbits. We also note that new islands are created inside the island centered
at the periodic points with period 6 when we increase $B$. Moreover, as $B$ continues growing, these new-formed islands get separated from the main island and form several isolated islands. By slightly increasing $B$, a similar pattern, the birth and separation of new islands, is observed in
the phase space of the billiard table.

Figure \ref{d-158} demonstrates a 6-periodic example of the remaining type of periodic point that may appear for lemon billiards. The 12 prominent peripheral islands corresponding to the periodic point and its mirror point in reverse time.

\begin{figure}[]
   \includegraphics[width=80mm]{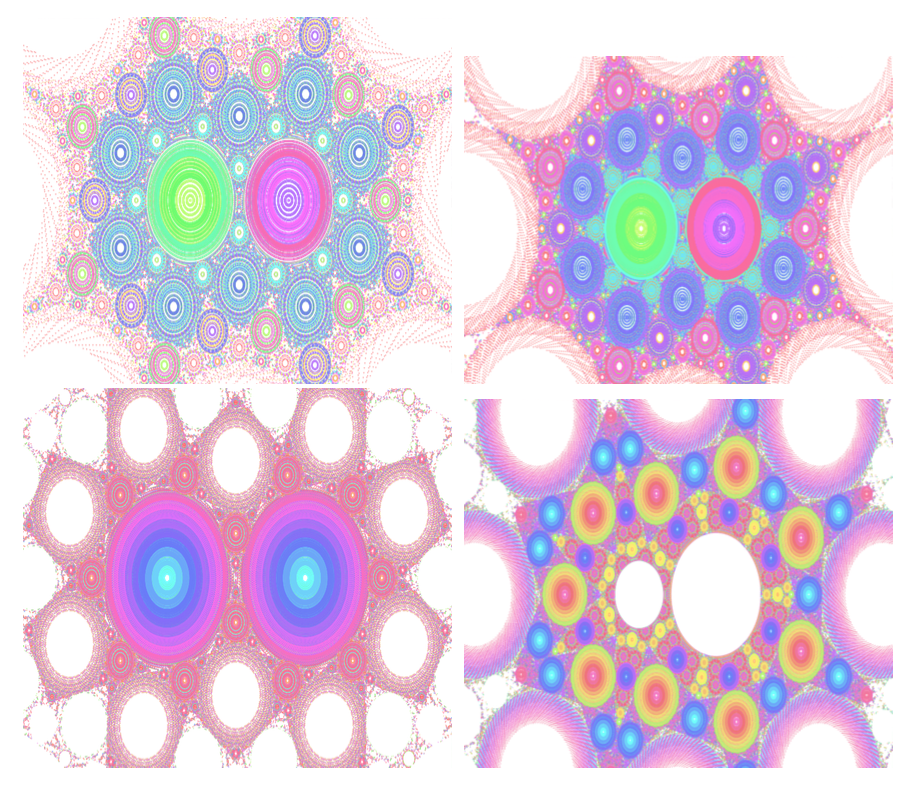}
   
 \caption{In umbrella billiards formed by small deformations of lemon billiards, the central 2-periodic point is replaced by two 2-periodic elliptic points surrounded by a multitude of higher order periodic points.  The scale of zooming on each phase portrait is proportional to the $B_1$ parameter: $B_1=0.01$ (upper left), $B_1=0.001$ (upper right), $B_1=0.00001$ (lower left), and $B_1=0.0000001$ (lower right), and accordingly only half of the elliptic islands are shown for the displayed points.}
 
\label{perturbation} 
 
\end{figure}

All of these types of periodic points may appear in umbrella billiards formed from the lemon modification. Specifically, if all of the collisions occur at boundary points removed from new corners then they will persist as $B_1$ increases. However, the dynamics may shift rapidly for umbrella billiards when the new corner appears near the periodic point, even for extremely small $B_1$ values. Figure \ref{perturbation} shows the details of the dynamical shift on several scales for the stable 2-periodic point shown in Figure \ref{d-135} when $B_1$ is small but nonzero. Notice that two elliptic 2-periodic points replace the original single 2-periodic point.

\subsection{Periodic points of moon type}

\begin{figure}[h!]
        \centering
\vspace{-.12in}        
\hspace*{-1.79cm}

        \includegraphics[width=80mm]{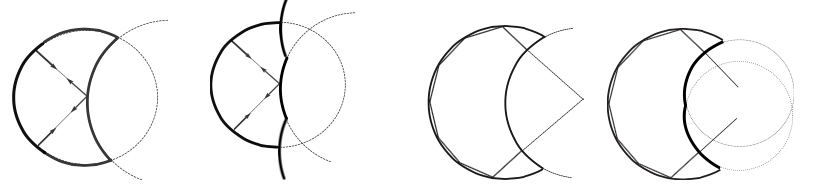}

 \caption{\small Periodic orbits for moon billiards containing radial trajectories may persist under the umbrella deformation.}
 \label{pers4e6}
 \end{figure}

In moon billiards, we have only one-sided periodic trajectories and no two-sided periodic orbits as observed in lemon billiards. Figure \ref{pers4e6} illustrates the existence of some elliptical points that experience several consecutive sliding collisions on the boundary
component of the unit disc, then collide perpendicularly on the boundary of the disc of radius $R$ and return after that. The requisite condition for this reversing collision is the extended trajectory passing through the center of one of the framing disks. This may also occur in billiards with the umbrella modification, and  these  periodic points will also occur.

\FloatBarrier

\vspace{-.15in}
\section{The transition to ergodicity in  umbrella billiards}
\label{sec:phase}
\begin{figure}[b]
  \centering
  
   \includegraphics[width=68mm]{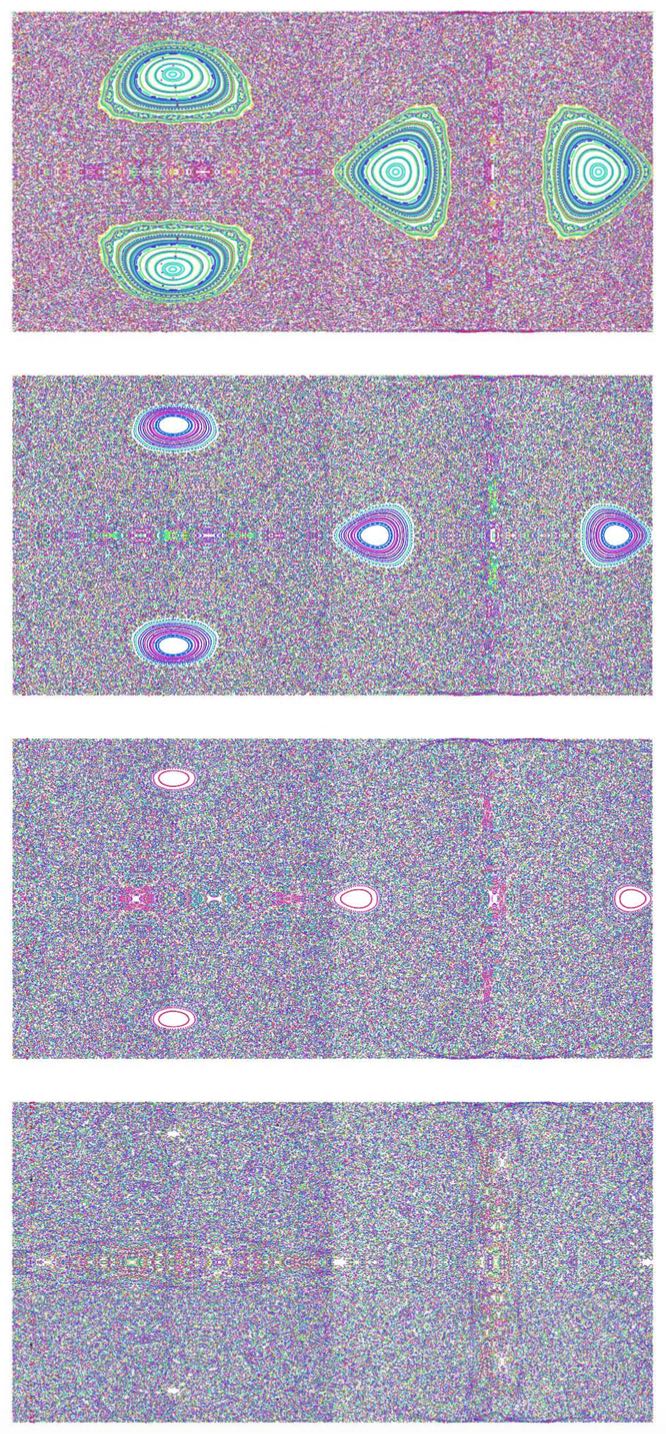}

\caption{Umbrella deformations of lemon type billiards $Q_L(1.27, 1.01, B_1)$ become increasingly chaotic as $B_1$ increases. Top: $0.4,$  middle top: $0.5,$ middle bottom: $0.6$, bottom: $0.7$.}\label{R127B101to08}
\label{lemtran}
\end{figure}

In this section we investigate the transition to ergodicity in the three-parameter umbrella families. 
For the two-parameter moon billiards investigated in  [\onlinecite{CZ15}], an ergodic boundary is hypothesized in parameter space demarcating a region of apparently ergodic billiards. (See the boundary in Figure \ref{moonoverview}.)  First, we will  look at examples $Q(R_0,B_0,0)$ on the non-ergodic side of the boundary and consider the effect of increasing $B_1$ by looking at the phase portraits. These examples have elliptic islands, the dissipation of which suggests that the corresponding umbrella billiards become ergodic.  Secondly, we consider the transition more broadly by mapping the ergodic border in the $\theta$ parameter space of moon type umbrella billiards for several values of $B_1$.

\begin{figure}[b]
  \includegraphics[width=68mm]{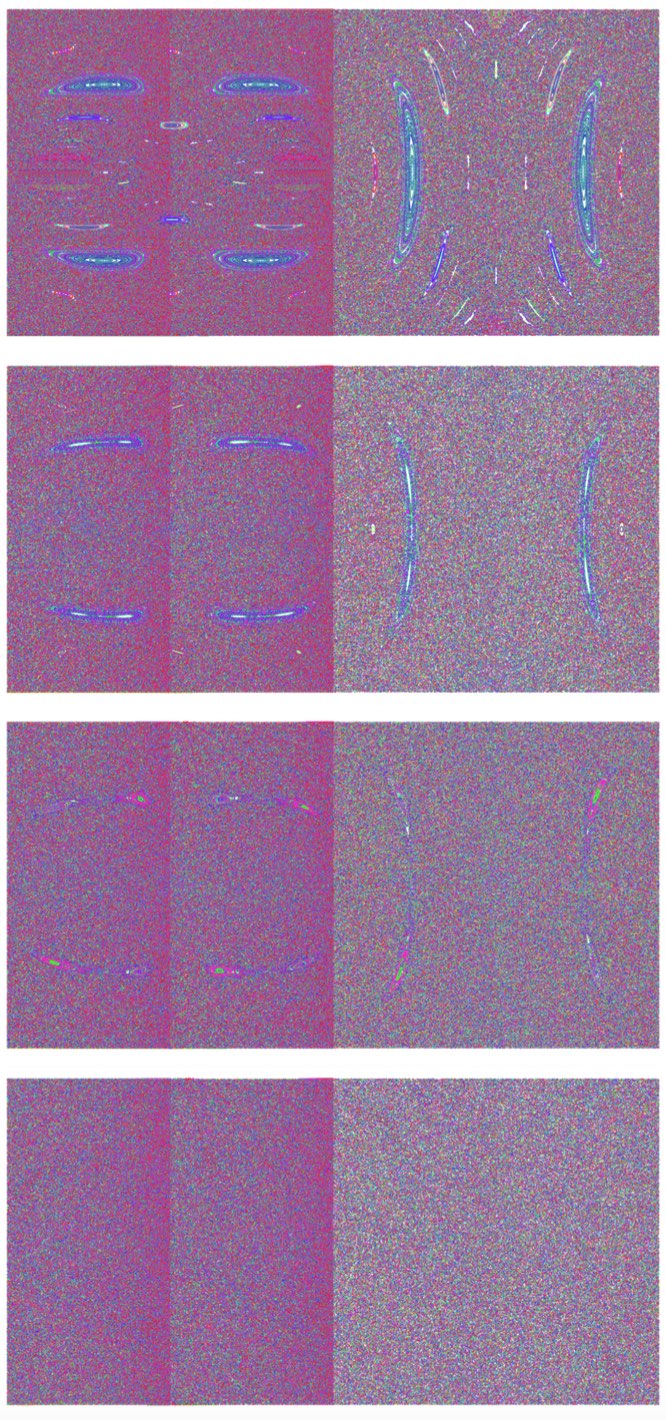}
     
 \caption{For umbrella deformations of type 1 moon billiards $Q_M^1(1.2,1.3,B_1),$ elliptic islands around periodic points vanish. Top: $B_1=0,$  middle top: $0.05,$ middle bottom: $0.1$, bottom: $0.2$.}
 
\label{moon1tran} 
\end{figure}

\begin{figure}[]
   \includegraphics[width=70mm]{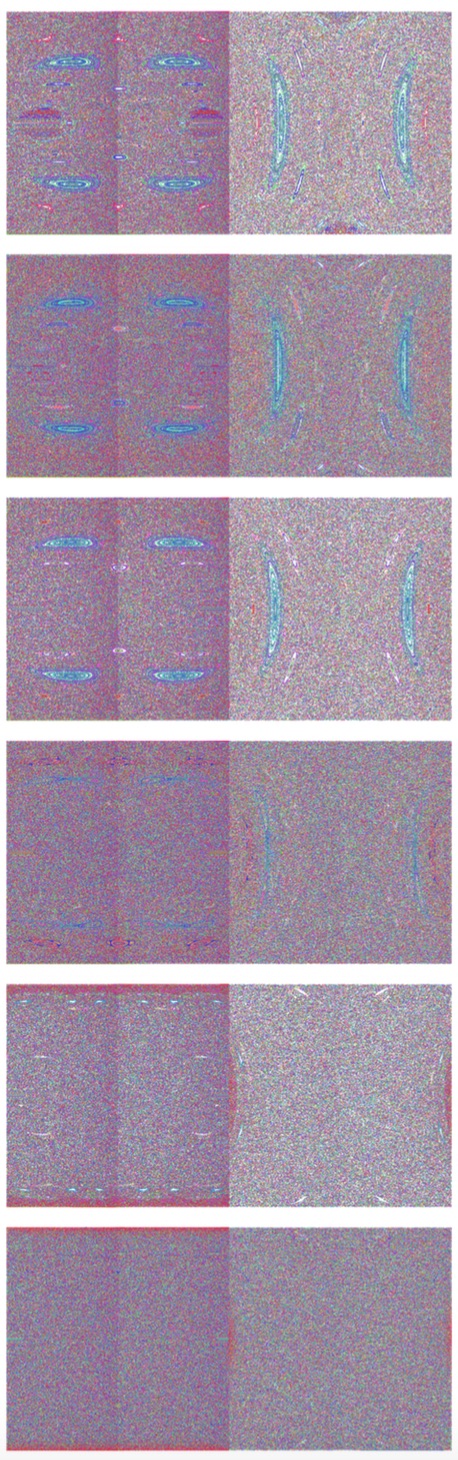}
\caption{$Q_M^2(1.2,1.3,B_1)$ with $B_1$ varying from $0$ to $1.8$. For type 2 moon based umbrella billiards, the transition towards chaotic behavior is less rapid. }
\label{moon2tran}
\end{figure}

Figure \ref{lemtran} shows phase portraits for lemon type umbrella billiards $Q_L(1.27,1.01, B_1)$. The islands corresponding to  the  prominent  4-periodic point shrink as $B_1$ is increased.  Figures \ref{moon1tran} and \ref{moon2tran} show similar examples for the first and second moon type umbrella billiards $Q_M(1.2,1.3,B_1)$, with a base moon type in a nonergodic region below the boundary hypothesized in [\onlinecite{CZ15}]. Both transition to a billiard that appears to be ergodic, but the transition occurs more rapidly for the first moon type.

\begin{figure}[b]
 \includegraphics[width=80mm]{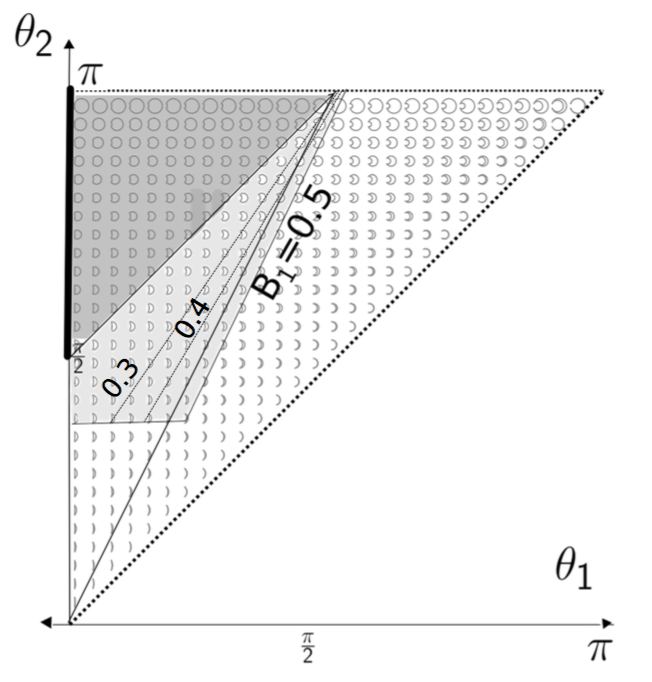}
\caption{The known region of numerical ergodicity for moon billiards (dark gray) is extended for the corresponding umbrella billiards (dark gray) when $B_1=0.5$. Smaller but observable corresponding regions exist for $B_1=0.3$ and $B_1=0.4$. }
\label{fig:moonumb}
\end{figure}

The greater efficacy of the first type in generating chaotic behavior is typical among the examples we have investigated, and accordingly it is a candidate for examining the shift in the ergodic boundary.  For $B_1=0$, the previously conjectured boundary translated to $\theta$ parameters is the line $\theta_2=\theta_1+\frac{\pi}{2}$. (See Figure \ref{2theta}.)  For a fixed $B_1>0$, we seek the boundary value for a fixed $\theta_2$ by sampling phase portraits of varying $\theta_1$ outside of the ergodic region for the base ($B_1=0$) case.  A transitional boundary can be identified by approximating the points $(\theta_1, \theta_2)$ near which the stable elliptic islands vanish and then adjusting $\theta_1$ to a smaller scale to refine the estimate. This approach cannot rule out the possibility of other tiny elliptic islands existing, but does give a candidate for a transitional boundary. Additionally, this technique was feasible in the selected region in part due to the fact that the persistent elliptic periodic points were few and of low order; in regions of parameter space where this is not the case the transition is less clear. Here, however, the boundary may be identified with a high level of precision. Specifically, each data point was obtained by isolating two phase portraits with a $\theta_1$ separation of $0.001$, one of which showed a clearly visible elliptic island while the other appeared to be ergodic with no evidence of elliptic islands.  

Using this method to estimate the boundary, the transition curves appear to be linear in $\theta$ for umbrella billiards of the first moon type. The values obtained as described fit a linear model to such a degree of accuracy that we conjecture the transition is linear. However, at this time it is not  clear why this should be true. 
Figure \ref{fig:moonumb} shows the boundaries for three values of $B_1$, giving evidence of an extended range of ergodic parameters compared to the $B=0$ plane corresponding to unmodified moon and lemon billiards.

\section{The Lyapunov exponents for $Q_M^1(1,B,B_1)$ and $Q_L(1,B,B_1)$}
\label{sec:lyap}

\begin{figure}[b]
\centering 
\includegraphics[width=85mm]{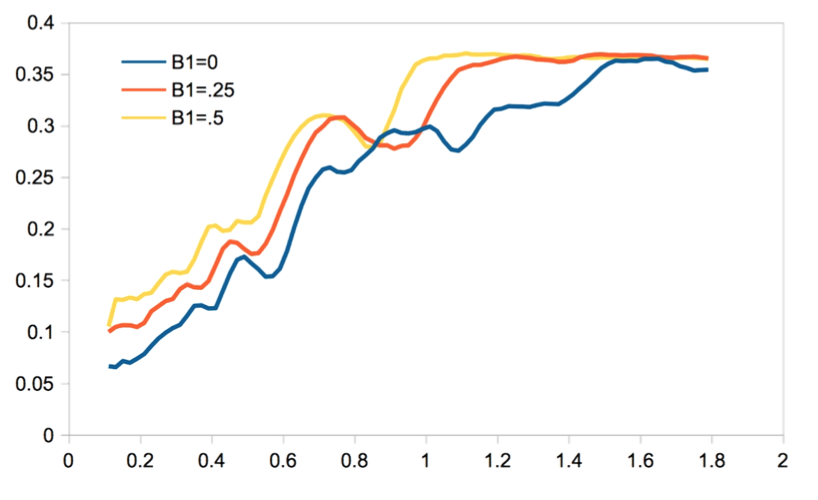}
\caption{\small Graphs of Lyapunov exponent $\overline{\lambda}$ as function of $B$ for three values of  $B_1$ in moon type 1 umbrella billiards. In all cases $R=1$. As $Q_M(1,B,0)$ already exhibits chaotic behavior for $B>\sqrt{2}$, the effect of increasing $B_1$ is not significant for higher value of $B$.}
\label{lyap1}
\end{figure}

\begin{figure}[t]
\includegraphics[width=85mm]{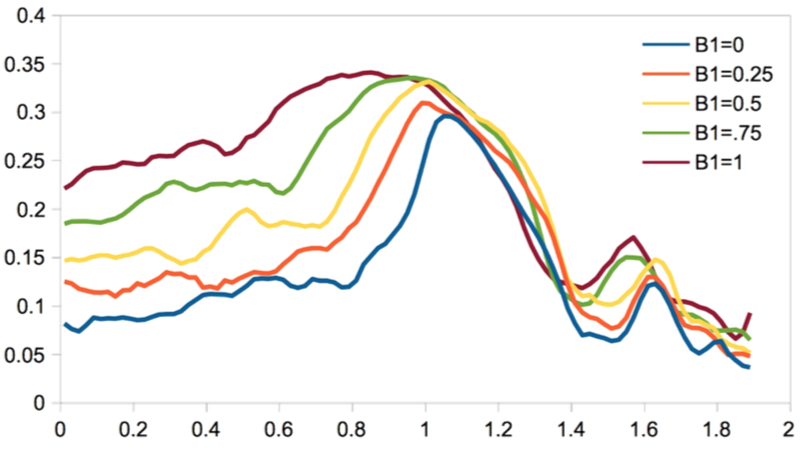}

\caption{\small Graphs of Lyapunov exponent $\overline{\lambda}$ as function of $B$ for  five values of $B_1$ for lemon type umbrella billiards. Notice that $Q_L(1,1,0)$ marks a transition at which the central elliptic periodic point becomes parabolic and the nature of the dynamics shifts. The effect of increasing $B_1$ appears to decrease after this transition.}
\label{lyap2}
\end{figure}

The observations in Section \ref{sec:phase} suggest that in many cases a transition towards chaotic behavior occurs as $B_1$ increases. In this section we wish to quantify the transition by calculating the Lyapunov exponents.  The Lyapunov exponent provides a quantitative measure of the stability of trajectories and has been widely used to quantify the average expansion or contraction rate for a small volume of initial conditions [\onlinecite{CM06}].  On one hand, if the maximum Lyapunov exponent is not positive, we have an indication of regularity and the dynamics can be periodic. On the other hand, if at least one Lyapunov exponent is positive, the orbit is said to be unstable and chaotic. Thus, introducing this measure of chaos we can compare tables for different parameters. As our primary interest is in the umbrella classes, we will consider values of $B_1$, looking at the Lyapunov exponents as $B$ varies with $R=1$ throughout this section.

The positive Lyapunov exponent, corresponding to the direction of expansion, is given by
$$\lambda^u=\underset{n\to\infty}{\lim_{\|\xi_0\|\to 0}}\;\frac{1}{n}\ln \frac{\|\eta_n\|}{\|\eta_0\|},$$
where $\eta_n$ is the $n$-distance in phase space between a fixed orbit and a nearby orbit  that begins with a nearby initial condition $\eta_0$. To numerically estimate,  one may use the finite-time Lyapunov indicator (LI):
$$\lambda(\eta_0,n)=\frac{1}{n}\ln \frac{\|\eta_n\|}{\|\eta_0\|}.$$
The basic features of the phase space of a chaotic system  can be discovered very quickly by calculating a large number of LIs for short time, and for parameters corresponding to ergodic billiard tables the Lyapunov exponents are constant almost everywhere and a random sampling for the initial conditions would suffice. However, for many parameters of lemon or moon type umbrella billiards large elliptic islands persist. To focus on the chaotic behavior, one might restrict choices of initial conditions to apparently ergodic regions of phase space, while a more systematic approach 
is given in [\onlinecite{SESF04}], systematically checking for quasiperiodic orbits with zero Lyapunov exponents and eliminating them from consideration. Since only we are interested in the relative change as $B_1$ varies, we use a third option. Define a scaled Lyapunov exponent by
$$\overline{\lambda}_{k,n} =\frac{1}{k^2}\sum \lambda(\eta_j,n),$$
where the initial conditions $\eta_j$ are uniformly sampled using a $k \times k$ grid of phase space. Hence, the values of $\overline{\lambda}$ will be lower than the true values of the Lyapunov exponents of the  chaotic region, but will still provide information about the relative change as $B_1$ increases. 

For $\| \Delta x \| = 10^{-6}$, $k=40$, and $n=10$, we calculate $\overline{\lambda}$ for $R=1$, $B$ varying across the trial, and $B_1$ varying between trials. The results (Figures \ref{lyap1} and \ref{lyap2}) suggest that in these cases the umbrella billiards can be more chaotic than their base moon and lemon types.

\section{Acknowledgments}

H-K. Zhang is partially supported by  NSF grant DMS-1151762. M.F. Correia is partially supported by the research center CIMA-UE, FCT Portugal funding program. M. F. Correia express her sincere gratitude to the Department of Mathematics and Statistics at University of Massachusetts Amherst for their support and hospitality, where this paper was written during her stay.

\end{document}